\input amstex
\input Amstex-document.sty

\pageno 457

\def\shorttitle{Non-zero Degree Maps between 3-Manifolds}
\def\shortauthor{S. C. Wang}

\topmatter
\title\nofrills{\boldHuge Non-zero Degree Maps between 3-Manifolds*}
\endtitle

\author \Large Shicheng Wang\dag \endauthor

\thanks *Supported by grants of MSTC and NSFC. \endthanks

\thanks \dag Department of Mathematics, Peking University, Beijing 100871,
China. E-mail:
\endthanks

\thanks\noindent swang\@sxx0.math.pku.edu.cn
\endthanks

\abstract\nofrills \centerline{\boldnormal Abstract}

\vskip 4.5mm

{\ninepoint First the title could be also understood as ``3-manifolds related by non-zero degree maps" or
``Degrees of maps between 3-manifolds" for some aspects in this survey talk.

The topology of surfaces was completely understood at the end of 19-th century, but maps between surfaces kept to
be an active topic in the 20-th century and many important results just appeared in the last 25 years. The
topology of 3-manifolds was well-understood only in the later 20-th century, and the topic of non-zero degree maps
between 3-manifolds becomes active only rather recently.

We will survey questions and results in the topic indicated by the title, present its relations to 3-manifold
topology and its applications to  problems in geometry group theory, fixed point theory and dynamics.

There are four aspects addressed: (1) Results concerning the existence and finiteness about the maps of non-zero
degree (in particular of degree one) between 3-manifolds and their suitable correspondence about epimorphisms on
knot groups and 3-manifold groups. (2) A measurement of the topological complexity on 3-manifolds and knots given
by ``degree one map partial order", and the interactions between the studies of non-zero degree map among
3-manifolds and of topology of 3-manifolds. (3) The standard forms of non-zero degree maps and automorphisms on
3-manifolds and applications to minimizing the fixed points in the isotopy class. (4) The uniqueness of the
covering degrees between 3-manifolds and the uniqueness  embedding indices (in particular the co-Hopfian property)
between Kleinian groups .

The methods used are varied, and we try to describe them briefly.

\vskip 4.5mm


}
\endabstract
\endtopmatter

\document

\baselineskip 4.5mm \parindent 8mm

\specialhead \noindent {\boldLARGE 0. Introduction} \endspecialhead

The topology of surfaces was completely understood by the end of
19th century, but maps between surfaces keep to be an active topic
in the 20th century, and some basic results just appeared in the
last 25 years, among which are the Nielsen-Thurston classification
of surface automorphisms [Th3], and Edmonds' standard form for
surface maps [E1]. Then fine results followed, say the realization
of Nielsen number in the isotopy class of surface automorphisms by
Jiang [J], and the simple loop theorem for surface maps by Gabai
[Ga2].

The topology of 3-manifolds was well-understood only in the later
20th century  due to many people's deep results,  in particular
Thurston's great contribution to the geometrization of
3-manifolds, and the topic of non-zero degree maps between
3-manifolds becomes active only rather recently.

We will survey the results and questions in the topic indicated by
the title, present its relations to 3-manifold topology and its
applications to problems in geometry group theory, fixed point
theory and dynamics. The methods used are varied, and we try to
describe them briefly.

For standard terminologies of 3-manifolds and knots, see the
famous books of J.\ Hempel, W.\ Jaco and D.\ Rolfsen. For a proper
map $f:M\to N$ between oriented compact 3-manifolds,
$\text{deg}(f)$, the degree of $f$, is defined in most books of
algebraic topology. A closed orientable 3-manifold is said to be
{\it geometric} if it admits one of the following geometries:
$H^3$ (hyperbolic), $\widetilde {PSL_2 R}$, $H^2\times E^1$, Sol,
Nil, $E^3$ (Euclidean), $S^2\times E^1$, and $S^3$ (spherical).  A
compact orientable 3-manifold $M$ admits a {\it geometric
decomposition\/} if each prime factor of $M$ is either geometric
or Haken. Thurston's geometrization conjecture asserts that any
closed orientable 3-manifold admits a geometric decomposition.
Each Haken manifold $M$ with $\partial M$ a (possibly empty) union
of tori has a Jaco-Shalen-Johannson (JSJ) torus decomposition,
that is, it contains a minimal set of tori, unique up to isotopy,
cutting $M$ into pieces such that each piece is either a Seifert
manifold or a simple manifold, which admits a complete hyperbolic
structure with finite volume [Th2].

In the remainder of the paper,  all manifolds are assumed to be
{\bf compact and orientable}, all automorphisms are {\bf
orientation preserving}, all knots are {\bf in $S^3$}, all
Kleinian groups are {\bf classical}, and all maps are {\bf
proper,} unless otherwise specified.

 Let $M$ and $N$ be 3-manifolds and $d>0$ an
integer.  We say that $M$ {\it $d$-dominates\/} (or simply
dominates) $N$ if there is a map $f: M\to N$ of degree $\pm d$.
Denote by $D(M,N)$ the set of all possible degrees of maps from
$M$ to $N$.    A 3-manifold $M$ is {\it small} if each closed
incompressible surface in $M$ is boundary parallel.

Due to space limitation, quoted literature are only partly listed
in the references; while the others are briefly indicated in the
context.

\specialhead \noindent {\boldLARGE 1. Existence and finiteness}
\endspecialhead

A fundamental question in this area (and in 3-manifold theory) is
the following.

\proclaim{Question 1.1} \rm  Given a pair of closed 3-manifolds
$M$ and $N$, can one decide if $M$ $d$-dominates $N$?  In
particular, can one decide if $M$ 1-dominates $N$?  \endproclaim

The following two natural problems concerning finiteness can be
considered as testing cases of Question 1.1.

\proclaim{Question 1.2 {\rm[Ki, Problem 3.100 (Y. Rong)]}}\rm Let
$M$ be a closed $3$-manifold.  Does $M$ 1-dominate at most
finitely many closed 3-manifolds?  \endproclaim

\proclaim {Question 1.3}\rm Let $N$ be a closed 3-manifold.  When
is $|D(M,N)|$ finite for any closed 3-manifold $M$?  \endproclaim

An important progress towards the solution of Question 1.2 is the
following:

\proclaim{Theorem 1.1 {\rm ([So2], [WZh2], [HWZ3])}} Any closed
$3$-manifold $1$-dominates at most finitely many geometric
$3$-manifolds.
\endproclaim

Theorem 1.1 was proved by Soma when the target manifolds $N$ admit
hyperbolic geometry [So2].  The proof is based on the argument of
Thurston's original approach to the deformation of acylindrical
manifolds.  Porti and Reznikov had a quick proof of Soma's result,
based on the volume of representations [Re2].  However Soma's
approach deserves attention as it proves that the topological
types of all hyperbolic pieces in closed Haken manifolds
1-dominated by $M$ are finite [So3].  Theorem 1.1 was proved in
[WZh2] when the target manifolds admit geometries of $H^2\times
E^1$, $\widetilde{PSL_2(R)}$, Sol or Nil.  The proof for the case
of $H^2\times E^1$ geometry invokes Gabai's result that embedded
Thurston Norm and singular Thurston Norm are equal [Ga1], and the
proof for case of $\widetilde {PSL_2(R)}$ geometry uses Brooks and
Goldman's work on Seifert Volume [BG].  Theorem 1.1 was proved in
[HWZ3] when the target manifolds admit $S^3$ geometry, using the
linking pairing of 3-manifolds.  Note that only finitely many
3-manifolds admit the remaining two geometries.

For maps between 3-manifolds which are not necessarily orientable,
there is a notion of geometric degree (See D. Epstein, Proc.
London Math. Soc. 1969).  It is worth mentioning that if
$d$-dominating maps are defined in terms of geometric degree, then
 Rong  constructs a non-orientable 3-manifold which 1-dominates
infinitely many lens spaces [Ro3]. Actually there is a
non-orientable hyperbolic 3-manifold which 1-dominates infinitely
many hyperbolic 3-manifolds [BW1].  Such examples do not exist in
dimension $n>3$ due to Gromov's work on simplicial volume and
H.C.\ Wang's theorem that, for any $V>0$, there are at most
finitely many closed hyperbolic $n$-manifolds of volume $<V$.

The answer to Question 1.2 is still unknown for closed irreducible
3-manifolds admitting geometric decomposition.  The following
result is related.

\proclaim{Theorem 1.2 {\rm ([Ro1], [So3])}} For any 3-manifold $M$
there exists an integer $N_{M}$, such that if $M = M_0\to M_1\to
....\to M_k$ is a sequence of degree one maps with $k>N_{M}$, and
each $M_i$ admits a geometric decomposition, then the sequence
contains a homotopy equivalence.  \endproclaim

The situation for Question 1.3 can be summarized in the following
theorem.

\proclaim{Theorem 1.3 {\rm ([Gr1], [BG], [W2])}} Suppose $N$ is a
closed 3-manifold admitting geometric decomposition. Then

(1) $|D(M,N)|$ is finite if either a prime factor of $N$ contains
a hyperbolic piece in its JSJ decomposition, or $N$ itself admits
the geometry of $\widetilde {PSL_2(R)}$.

(2) $|D(N,N)|$ is infinite if and only if either (i) $N$ is
covered by a torus bundle over the circle or a surface$\times
S^1$, or (ii) each prime factor of $N$ has a cyclic or finite
fundamental group.
\endproclaim

Part (1) of Theorem 1.3 follows from the work of Gromov [Gr1] and
Brooks-Goldman [BW]. Part (2) can be found in [W2].  Note that if
$|D(N, N)|$ is infinite and $D(M,N)$ contains non-zero integers,
then $|D(M, N)|$ is also infinite. I suspect that Theorem 1.3 (2)
indicates a general solution to Question 1.3.

There are many partial results for Question 1.1: When both $M$ and
$N$ are Seifert manifolds with infinite fundamental groups Rong
has an algorithm to determine if $M$ 1-dominates $N$ [Ro3]. When
$N$ is the Poincare homology sphere and a Heegaard diagram of $M$
is given, Hayat-Legrand, Matveev and Zieschang have an algorithm
to decide if $M$ $d$-dominates $N$ [HMZ]. There are simple answers
to Question 1.1 in the following cases: (1) $M$ and $N$ are prism
spaces and $d=1$ [HWZ2]; (2) $M=N$ admit geometry of $S^3$ and
$f_*$ an automorphism on $\pi_1$ [HKWZ]; (3) $N$ is a lens space.
I will state (3) as a theorem, since both its statement and proof
are short, and since it has rich connections with previous results
and with different topics.

\proclaim{Theorem 1.4 {\rm ([HWZ1], [HWZ3])}} A closed
$3$-manifold $M$ $d$-dominates the lens space $L(p,q)$ if and only
if there is an element $\alpha$ in the torsion part of $H_1(M,
{\Bbb Z})$ such that $\alpha \odot \alpha =\frac {dq}p$ in ${\Bbb
Q}/{\Bbb Z}$, where $\alpha \odot \alpha$ is the self-linking
number of $\alpha$.
\endproclaim

A direct consequence of Theorem 1.4 is the known fact that
$L(p,q)$ 1-dominates $L(m,n)$ if and only if $p=km$ and $n=kqc^2$
mod $m$. This fact has at least four different proofs: using
equivariant maps between spheres by de Rham (J. Math. 1931) and by
Olum (Ann. of Math. 1953), using Whitehead torsion by Cohen (GTM
10,  1972), using  pinch in [RoW] and using linking pair in
[HWZ1].

Degree one maps from general 3-manifolds to some lens spaces, in
particular  the $RP^3$, have been studied by Bredon-Wood (Invent.
Math. 1969) and by Rubinstein (Pacific J. Math. 1976) to find
one-sided incompressible surfaces, by Luft-Sjerve (Topo. Appl.
1990) to study cyclic group actions on 3-manifolds, by
Shastri-Williams-Zvengrowski [SWZ] in theoretical physics, by
Taylor (Topo. Appl. 1984) to define normal bordism classes of
degree one maps, and by Kirby-Melvin (Invent. Math. 1991) to
connect with new 3-manifold invariants.

Degree one maps induce epimorphisms on $\pi_1$. There are easy
examples  indicate that Question 1.2 does not have direct
correspondence in the level of 3-manifold groups [BW1], [RWZh].
However the following related question was raised in 1970's.

\proclaim {Question 1.4 {\rm [Ki, Problem 1.12 (J. Simon)]}}\rm
Conjectures:

(1) Given a knot group $G$, there is a number $N_G$ such that any
sequence of epimorphisms of knot groups $G\to G_{1}\to .... \to G_{n}$
with $n\ge N_G$ contains an isomorphism.

(2) Given a knot group $G$, there are only finitely many knot groups
$H$ for which there is an epimorphism $G\to H$.
\endproclaim

According to a conversation with Gonzalez-Acuna, who discussed
Question 1.4 with Simon before it was posed, the epimorphisms in
Question 1.4 are peripheral preserving in their minds.

\proclaim{Theorem 1.5 {\rm ([So5], [RW])}} The conjecture in
Question 1.4 (1) holds if the knot complements involved are small.
The conjecture in Question 1.4 (2) holds if the knot complements
are small and the epimorphisms are peripheral preserving.
\endproclaim

The first claim is due to Soma [So5] and the second claim  is in
[RW]. Both of them invoke Culler-Shalen's work on the
representation varieties of knot groups. It is also proved that
any infinite sequence of epimorphisms among 3-manifold groups
contains an isomorphism if all manifolds are either hyperbolic
[So5] or Seifert fibered [RWZh]. In [RWZh], the proof uses the
fact that epimorphisms between aspherical Seifert manifolds with
the same $\pi_1$ rank are realized by maps of non-zero degree.
Both this fact and Question 1.4 (1) are variations of the Hopfian
property.

We end this section by mention that there are results about
$D(M,N)$ in [DW] for $(n-1)$-connected $2n$-manifolds, $n>1$,
which are quite explicit  and of interest from both topological
and number-theoretic point of view.

\specialhead \noindent {\boldLARGE 2. Uniqueness }
\endspecialhead

The following question is raised in 1970's.

\proclaim{Question 2.1 {\rm [Ki, Problem 3.16 (W. Thurston)]}} \rm
Suppose a 3-manifold $M$ is not covered by (surface)$\times S^1$
or a torus bundle over $S^1$. Let $f,g :M\to N$ be two coverings,
must $deg(f)=deg(g)$?  \endproclaim

It is known [WWu2] that Question 2.1 has positive answer if $M$
admits geometric decomposition and is not a graph manifold ($M$ is
a graph manifold if each piece of its JSJ decomposition is Seifert
fibered.)  For graph manifolds there are four different covering
invariants  introduced in middle 1990's by [WWu2],  Luecke and Wu
[LWu],  Neumann [N] and   Reznikov [Re1] . Unfortunately all those
four covering invariants are either vanishing or not well-defined
for some non-trivial graph manifolds. It is also known that
covering degree is uniquely determined if the graph manifold in
the target is either a knot complement [LWu] or its corresponding
graph is simple [WWu2, N]. The positive answer to Question 2.1 for
graph manifolds was finally obtained in [YW], using the matrix
invariant defined in [WWu2] and an elegant application of matrix
theory due to Yu.

\proclaim{Theorem 2.1 {\rm ([WWu2], [YW])}} For 3-manifolds
admitting geometric decomposition and not covered by either
(surface)$\times S^1$ or a torus bundle over $S^1$, covering
degrees are uniquely determined by the manifolds involved.
\endproclaim

It is worth mentioning an interesting fact that any knot
complement is non-trivially covered by at most two knot
complements and any knot complement non-trivially covers at most
one knot complement. The first claim follows from the cyclic
surgery theorem of Culler-Gordon-Luecke-Shalen and the positive
answer to the Smith Conjecture. The second claim is in [WWu1].

Question 2.1 is equivalent to asking the uniqueness of indices of
finite index embeddings between 3-manifold groups. Recently there
are also some discussions on the uniqueness of indices of
self-embeddings of groups.  A group $G$ is said to be co-Hopf if
each self-monomorphism of $G$ is an isomorphism.

\proclaim{Question 2.2}\rm Let $G$ be either a 3-manifold group,
or a Kleinian group, or a word hyperbolic group.  When is $G$
co-Hopf?
\endproclaim

The cohopficity of groups were first considered by Baer (Bull. AMS
1944).  For word hyperbolic groups it was first considered by
Gromov in 1987 [Gr2, p.157], and subsequently by Rips-Sela (GAFA,
1994), Sela [Se], and Kapovich-Wise (Isreal J. Math. 2001).
Cohopficity of 3-manifold groups was first studied in 1989 by
Gonzalez-Acuna and Whitten [GWh], and then in [WWu2] and [PW]. The
answer for 3-manifolds admitting geometric decomposition with
boundary either empty set or a union of tori is known [GWh],
[WWu2], and partial results for 3-manifolds with boundary of high
genus surfaces are in [PW]. Cohopficity of Kleinian groups was
first considered in 1992 in an early version of [PW], then in 1994
in an early version of [WZh1], also by Ohshika-Potyagailo (Ann.
Sci. Ecole Norm. Sup. 1998) and Delzant-Potyagailo (MPI Preprint,
2000) for high dimensional Kleinian groups.

\proclaim{Theorem 2.2} Suppose $K$ is a non-elementary, freely
indecomposable, geometrically finite Kleinian group and $K$
contains no  $Z\oplus Z$ subgroup. Then

(1){\rm [Se], [PW], [WZh1]} $K$ is co-Hopf if $K$ is a group of
one end.

(2){\rm [WZh1]} If the singular locus of the hyperbolic 3-orbifold
$H^3/K$ is a 1-manifold, then $K$ is co-Hopf if and only if no
circle component of singular locus meets a minimal splitting
system of hyperbolic cone planes.  \endproclaim

The proof of Theorem 2.2 (1) in [WZh1], influenced by that of
torsion free case in [PW],  use a generalization of
Thurston-Gromov's finiteness theorem on the conjugacy classes of
group embeddings (Delzant, Comm. Math. Helv. 1995) and a proper
conjugation theorem of Kleinian groups (Wang-Zhou, Geometriae
Dedicata, 1995). Theorem 2.2 (2) is proved by using 3-dimensional
hyperbolic orbifold structures and orbifold maps, which turn out
to be useful geometric tools.

Note that groups in Theorem 2.2 are word hyperbolic groups.
According to Sela ([Se] and his MSRI preprints in 1994), people
once expected that a non-elementary word hyperbolic group is
co-Hopf if and only if it has one end. Sela proved this
expectation for the torsion free case [Se]. Theorem 2.2 (2) and
examples in [WZh1] show that cohopficity phenomenon is very
complicated in the torsion case. In particular there are co-Hopf
word hyperbolic groups which have infinitely many ends.

Inspired by Questions 2.1 and 2.2 it is natural to ask

\proclaim{Question 2.3}\rm  Are the indices (including the
infinity) of embeddings $H\to G$ between co-Hopf groups unique?
\endproclaim

\specialhead \noindent {\boldLARGE 3. Interactions with 3-manifold topology}\endspecialhead

Degree one maps define a partial order on Haken manifolds and
hyperbolic 3-manifolds. By Gordon-Luecke's theorem knots are
determined by their complements [GL].  We say that a knot $K$
1-dominates a knot $K'$ if the complement of $K$ 1-dominates the
complement of $K'$. 1-domination among knots also gives a partial
order on knots. This partial order seems to provide a good
measurement of complexity of 3-manifolds and knots. The reactions
of non-zero degree maps between 3-manifolds and 3-manifold
topology are reflected in the following very flexible

\proclaim{Question 3.1} \rm Suppose $M$ and $N$ are 3-manifolds
(knots) and $M$ 1-dominates ($d$-dominates) $N$.

(1) Is $\sigma(M)$ not ``smaller" than $\sigma(N)$ for a
topological invariant $\sigma$?

(2) If $M$ and $N$ are quite ``close", are they homeomorphic?  do
they admit the same topological structure?  \endproclaim

Positive answers to Question 3.1 (1) are known in many cases.
Suppose $M$ 1-dominates $N$. Then $\sigma(M)\ge \sigma (N)$ when
$\sigma$ is either the rank of $\pi_1$, or Gromov's simplicial
volume, or Haken number (of incompressible surfaces), or genus of
knots; $\sigma (N)$ is a direct summand of $\sigma(M)$ when
$\sigma$ is the homology group, and $\sigma(N)$ is a factor of
$\sigma(M)$ if $\sigma$ is the Alexander polynomial of knots. The
answer to Question 3.1 (1) is still unknown for many invariants of
knots and 3-manifolds, for example crossing number, unknotting
number, Jones polynomial, knot energy, and tunnel number, etc. Li
and Rubinstein are specially interested in Question 3.1 (1) for
Casson invariant in order to prove   it is a homotopy invariant
[LRu].

There are both positive and negative answers to Question 3.1 (2),
depending on the interpretation of the problem.  On the negative
side, Kawauchi has constructed, using the imitation method
invented by himself, degree one maps between non-homeomorphic
3-manifolds $M$ and $N$ with many topological invariants
identical, see his survey paper [Ka]. On the positive side, there
are many results. An easy one is that if $M$ $d$-dominates $N$ and
both $M$ and $N$ are aspherical Seifert manifolds, then the Euler
number of $M$ is zero if and only if that of $N$ is zero [W1]. A
deeper result is Gromov-Thurston's Rigidity theorem, which says
that a degree one map between hyperbolic 3-manifolds of the same
volume is homotopic to an isometry [Th2]. The following are some
recent results in this direction.

\proclaim{Theorem 3.1 {\rm ([So4], [So1])}} (1)  For any $V>0$,
suppose $f:M\to N$ is a degree one map between closed hyperbolic
3-manifolds with $Vol (M)<V$. Then there is a constant $c=c(V)$
such that $(1-c) {\text Vol} (M)\le {\text Vol} (N)$ implies that
$f$ is homotopic to an isometry.

(2)   If $M\to N$ is a map of degree $d$ between Haken manifolds
such that $||M||=d||N||$, then $f$ can be homotoped to send $H(M)$
to $H(N)$ by a covering, where $||*||$ is the Gromov norm and
$H(*)$ is the hyperbolic part under the JSJ decomposition.
\endproclaim

\proclaim{Theorem 3.2 {\rm ([BW1], [BW2])}} (1)  Let $M$ and $N$
be two closed irreducible 3-manifolds with the same first Betti
number and suppose $M$ is a surface bundle. If $f:M\to N$ is a map
of degree $d$, then $N$ is also a surface bundle.

(2)  Let $M$ and $N$ be two closed, small hyperbolic 3-manifolds.
If there is a degree one map $f: M\to N$ which is a homeomorphism
outside a submanifold $H\subset N$ of genus smaller than that of
$N$, then $M$ and $N$ are homeomorphic.
\endproclaim

(1) and (2) of Theorem 3.1 provide a stronger version and a
generalization of Gromov-Thurston's Rigidity theorem,
respectively. In respect of Theorems 3.2, the following examples
should be mentioned: There are degree one maps between two
non-homeomorphic hyperbolic surface bundles with the same first
Betti number and between two non-homeomorphic small hyperbolic
3-manifolds [BW2]. The constructions of those maps are quite
non-trivial. There are many applications of Theorem 3.2. We list
two of them which are applications of Theorem 3.2 to Thurston's
surface bundle conjecture and to Dehn surgery respectively, where
degree one maps constructed by surgery on null-homotopic knots are
involved.

\proclaim{Theorem 3.3 {\rm ([BW1], [BW2])}} (1)  There are closed
hyperbolic 3-manifolds $M$ such that any tower of abelian covering
of $M$ contains no surface bundle.

(2)  Suppose $M$ is a small hyperbolic 3-manifold and that
$k\subset M$ is a null-homotopic knot, which is not in a 3-ball.
If the unknotting number of $k$ is smaller than the Heegaard genus
of $M$, then every closed 3-manifold obtained by a non-trivial
Dehn surgery on $k$ contains an incompressible surface.
\endproclaim

\specialhead \noindent {\boldLARGE 4. Standard forms }
\endspecialhead

\proclaim{Question 4.1} \rm What are standard forms of non-zero
degree maps and of automorphisms of 3-manifolds?  \endproclaim

Sample answers to analogs of Question 4.1 in dimension 2 are that
each map of non-zero degree between closed surfaces is homotopic
to a pinch followed by a branched covering [E1], and each
automorphism on surfaces can be isotoped to a map which is either
pseudo Anosov (Anosov), or periodic, or reducible [Th3].

\proclaim{Theorem 4.1 {\rm (Haken, Waldhausen, [E2], [Ro2])}} (1)
A degree one map between closed  3-manifolds is homotopic to a
pinch.

(2) A map of degree at least three between closed
3-manifolds is homotopic to a branched covering.

(3) A non-zero degree map between Seifert manifolds with infinite
$\pi_1$ is homotopic to a fiber preserving pinch followed by a fiber
preserving branched covering.  \endproclaim

(1) is proved by Haken (Illinois J. Math. 1966), also by
Waldhausen, and a quick proof using differential topology is in
[RoW]. (2) is proved by Edmonds [E2] quickly after
Hilden-Montesinos's result that each 3-manifold is a 3-fold
branched covering of 3-sphere. (3) is due to Rong [Ro2], which
invokes [E1].  According to conversations with D. Gabai and with
M. Freedman, people are still wondering if each map of degree 2
between closed 3-manifolds is homotopic to a pinch followed by a
double branched covering.

For non-prime 3-manifolds, Cesar de Sa and Rourke  claim that
every automorphism is a composition of those preserving and
permuting prime factors (Bull. AMS, 1979), and those so-called
sliding maps. A proof is given by Hendricks and Laudenbach [HL],
and by McCullough [Mc].

Standard forms of automorphisms on prime 3-manifolds admitting
geometric decomposition have been studied in [JWW]. The orbifold
version of Nielsen-Thurston's classification of surface
automorphisms is established, i.e., each orbifold automorphism is
orbifold-isotopic to a map which is either (pseudo) Anosov, or
periodic, or reducible.  We then have the following theorem.

\proclaim{Theorem 4.2 {\rm ([JWW])}} Let $M$ be a closed prime
3-manifold admitting geometric decomposition.  Let $f: M\to M$ be
an automorphism.  Let $\Cal T$ be the product neighborhood of the
JSJ tori. Then

(1) $f$ is isotopic to an affine map if $M$ is a 3-torus.

(2) $f$ is isotopic to an isometry if $M$ is the Euclidean
manifold having a Seifert fibration over $RP^2$ with two singular
points of index 2.

(3) $f$ is isotopic to a map which preserves the torus bundle
structure over 1-orbifold if $M$ admits the geometry of Sol.

(4) for all the remaining cases, $f$ can be isotoped so that $\Cal
T$ is invariant under $f$, and for each $f$-orbit $O$ of the
components in $\{\Cal T, \overline {M-\Cal T}\}$, $f|O$ is an
isometry if $O$ is hyperbolic, $f|O$ is affine if $O$ belongs to
$\Cal T$, otherwise there is a Seifert fibration on $O$ so that
$f$ is fiber preserving and the induced map on the orbifold is
either periodic, or (pseudo) Anosov, or reducible.
\endproclaim

As in dimension 2, standard forms in Theorems 4.2 are useful in
the study of fixed point theory and dynamics of 3-manifold
automorphisms.  The following is a result in this direction.

\proclaim{Theorem 4.3 {\rm ([JWW])}} Suppose $M$ is a closed
prime 3-manifold admitting geometric decomposition and $f: M\to M$
is an automorphism.  Then

(1) the Nielsen number $N(f)$ is realized in the isotopy class of $f$.

(2) $f$ is isotopic to a fixed point free automorphism unless some
component of the JSJ decomposition of $M$ is a Seifert manifold whose
orbifold is neither a 2-sphere with a total of at most three holes or
cone points nor a projective plane with a total of at most two holes
or cone points.
\endproclaim

\specialhead \noindent \boldLARGE References \endspecialhead

\widestnumber\key{HKWZ}

\ref \key BW1 \by M. Boileau and S.C. Wang \paper Non-Zero degree
maps and surface bundles over $S^1$ \jour J. Diff. Geom.  \vol 43
\pages 789--908 \yr 1996 \endref

\ref \key BW2 \by M. Boileau and S.C.Wang \paper Degree one maps,
incompressible surfaces and Heegaard genus \jour Preprint \yr 2002
\endref

\ref \key BG \by R. Brooks and W. Goldman \paper Volume in Seifert
space \jour Duke Math. J. \vol 51 \pages 529--545 \yr 1984 \endref

\ref \key CS \by M. Culler and P. Shalen \paper Varieties of group
representations and splittings of 3-manifolds \jour Ann. of Math.
\vol 117 \pages 109--146 \yr 1983 \endref

\ref \key DW \by H. Duan and S.C.Wang \paper The degree of maps
between manifolds \jour Math. Zeit. \toappear
\endref

\ref \key E1 \by A. Edmonds \paper Deformation of maps to branched
covering in dimension 2 \jour Ann. of Math. \vol 110 \pages
113--125 \yr 1979 \endref

\ref \key E2 \bysame \paper Deformation of maps to branched
covering in dimension 3 \jour Math. Ann. \vol 245 \pages 273--279
\yr 1979 \endref

\ref \key Ga1 \by D. Gabai \paper Foliations and the topology of
3-manifolds \jour J. Diff. Geom.  \vol 18 \yr 1983 \pages 479--536
\endref

\ref \key Ga2 \bysame \paper  Simple loop theorem \jour J. Diff.
Geom. \vol 21 \pages 143--149 \yr 1985 \endref

\ref \key GL \by C. Gordon and J. Luecke \paper Knots are
determined by their complements \jour JAMS \vol 2 \pages 371--415
\yr 1989 \endref

\ref \key GWh \by F. Gonzalez-Acuna and W. Whitten \paper
Embeddings of 3-manifold groups \jour  Mem. AMS \vol 474
\endref

\ref \key Gr1 \by M. Gromov \paper Volume and bounded cohomology
\jour Publ. Math. IHES  \vol 56 \pages 5--99 \yr 1983 \endref

\ref \key Gr2 \bysame \paper Hyperbolic groups \inbook Essays in
Group Theory \bookinfo edited by S.M. Gersten, MSRI Pub.  \vol 8
\publ Springer-Verlag, 75--263  \endref

\ref \key HMZ \by C. Hayat-Legrand, S. Matveev and H. Zieschang
\paper Computer calculation of the degree of maps into the
Poincare homology sphere  \jour Experiment. Math. \vol 10 \pages
497--508 \yr 2001
\endref

\ref \key HKWZ \by  C. Hayat-Legrand, E. Kudryavtseva, S.C.Wang
and H. Zieschang  \paper  Degrees of self-mappings of Seifert
manifolds with finite $\pi_1$ \jour  Rend. Istit. Mat. Univ.
Trieste \vol 32 \pages 131--147 \yr 2001
\endref

\ref \key HWZ1 \by C. Hayat-Legrand, S.C.Wang and H. Zieschang
\paper Degree one map onto lens spaces \jour Pacific J. Math. \vol
176 \pages 19--32 \yr 1996 \endref

\ref \key HWZ2 \bysame  \paper Minimal Seifert manifolds \jour
Math. Ann. \vol 308 \pages 673--700 \yr 1997 \endref

\ref \key HWZ3 \bysame \paper Any $3$-manifold $1$-dominates only
finitely many 3-manifolds supporting $S^3$ geometry \jour Proc.
AMS \vol 130 \pages 3117--3123 \yr 2002 \endref

\ref \key HL \by  H. Hendriks and F. Laudenbach \paper
Diffeomorphismes des sommes connexes en dimension trois. \jour
Topology \vol 23 \pages 423--443 \yr 1984
\endref

\ref \key J \by  B. Jiang \paper Fixed points of surface
homeomorphisms \jour Bull. AMS \vol 5 \pages 176--178 \yr 1981
\endref

\ref \key JWW \by B. Jiang, S.C.Wang, and Y-Q.Wu \paper
Homeomorphisms of 3-manifolds and the realization of Nielsen
number \jour Comm. Anal.  Geom. \vol 9 \pages 825--877 \yr 2001
\endref

\ref \key Ka \by A. Kawauchi \paper Topological imitations \inbook
Lectures at Knots 96 (ed. Shin'ichi Suzuki) \yr 1997, 19--37 \publ
World Sci. Publ. Co. \endref

\ref \key Ki \by  R. Kirby \paper Problems in low-dimensional
topology, \inbook Geometric topology, (ed. H.Kazez), AMS/ \publ
International Press  \yr 1997, 35--473 \endref

\ref \key LRu \by W. Li and J.H. Rubinstein \paper Casson
invariant is a homotopy invariant \jour Preprint \yr 2002
\endref

\ref \key LWu \by  J. Luecke and Y-Q. Wu \paper Relative Euler
number and finite covers of graph manifolds \inbook Geometric
topology, (ed. H.Kazez), AMS/ \publ International Press \vol 1
 \yr 1997, 80--103
\endref

\ref \key Mc \by  D. McCullough \paper Mapping of reducible
3-manifolds \inbook Proc. of the Semester in Geometric and
Algebraic Topology, Warsaw  \publ Banach Center, Publ. \yr 1986,
61--76
\endref

\ref \key N \by W.D. Neumann \paper Commensurability and virtual
fibration for graph manifolds \jour Topology \vol 36 \pages
355--378 \yr 1997 \endref

\ref \key PW \by  L. Potyagailo, and S.C. Wang \paper 3-manifolds
with co-Hopf fundamental groups \jour St. Petersburg Math. J.
(English Translation) \vol 11 \pages 861--881 \yr 2000 \endref

\ref \key Re1 \by  A. Reznikov \paper Volume of discrete groups
and topological complexity of homology spheres \jour Math. Ann.
\vol 306 \pages 547--554 \yr 1996 \endref

\ref \key Re2 \bysame \paper Analytic Topology \inbook Progress in
Math. Vol.201, 519--532. Birkhauser 2002
\endref

\ref \key RW \by  A. Reid and S.C.Wang \paper Non-Haken
3-manifolds are not large with respect to mappings of non-zero
degree \jour Comm. Anal. Geom. \vol 7 \pages 105--132  \yr 1999
\endref

\ref \key RWZh \by  A.Reid, S.C.Wang, Q.Zhou \paper Generalized
Hopfian property, a minimal Haken manifold, and epimorphisms
between 3-manifold groups\jour Acta Math. Sinica  \vol 18 \pages
157--172 \yr 2002 \endref

\ref \key Ro1 \by  Y. Rong \paper Degree one maps between
geometric 3-manifolds \jour Trans. AMS \yr 1992 \endref

\ref \key Ro2 \bysame \paper Maps between Seifert fibered spaces of infinite
$\pi_1$ \jour Pacific J. Math. \vol 160 \yr 1993 \pages 143--154 \endref

\ref \key Ro3 \bysame  \paper Degree one maps of Seifert manifolds and a
note on Seifert volume \jour  Topology Appl. \vol 64 \yr 1995
\pages 191--200 \endref

\ref \key RoW \by Y. Rong and S.C.Wang \paper The preimage of
submanifolds \jour Math. Proc. Camb.  Phil. Soc. \vol 112\pages
271--279 \yr 1992 \endref

\ref \key Se1 \by Z. Sela \paper Structure and rigidity in
(Gromov) hyperbolic groups \jour GAFA \vol 7 \pages 561--593 \yr
1997 \endref

\ref \key SWZ \by  A. Shastri, J.G.Williams and P.Zvengrowski
\paper Kinks in general relativity \jour Internat. J. Theor. Phys.
\vol 19 \pages 1--23 \yr 1980 \endref

\ref \key So1 \by T. Soma \paper A rigidity theorem for Haken
manifolds \jour Math. Proc.  Camb. Phil. Soc. \vol 118 \pages
141--160 \yr 1995 \endref

\ref \key So2 \bysame  \paper Non-zero degree maps onto hyperbolic
$3$-manifolds \jour J. Diff. Geom. \vol 49  \yr 1998 \pages
517--546 \endref

\ref \key So3 \bysame  \paper Sequence of degree-one maps between
geometric 3-manifolds. \jour Math. Ann. \vol 49  \yr 2000 \pages
733--742 \endref

\ref \key So4 \bysame \paper Degree one maps between hyperbolic
3-manifolds with the same limit \jour  Trans. AMS  \vol 353 \pages
2753  \yr 2001 \endref

\ref \key So5 \bysame \paper  Epimorphisms sequences between
hyperbolic 3-manifold groups \jour  Proc. AMS \vol 130 \yr 2002
\pages 1221--1223 \endref

\ref \key Th1 \by W. Thurston \book The Geometry and Topology of
Three-Manifolds \publ Princeton Lecture Notes \endref

\ref \key Th2 \bysame \paper Three dimensional manifolds, Kleinian
groups and hyperbolic geometry \jour Bull. AMS \vol 6 \pages
357--388 \yr 1982 \endref

\ref \key Th3 \bysame \paper On the geometry and dynamics of
diffeomorphisms of surfaces \jour Bull. AMS\vol 19 \pages 417--431
\yr 1988 \endref

\ref \key W1 \by S.C.Wang \paper The existence of non-zero degree
maps between aspherical 3-manifolds \jour Math. Zeit. \vol 208
\pages 147--160 \yr 1991 \endref

\ref \key W2 \bysame \paper The $\pi_1$-injectivity of self-maps
of non-zero degree on 3-manifolds \jour Math. Ann. \vol 297 \pages
171--189 \yr 1993 \endref

\ref \key WWu1 \by S.C.Wang and Y.Q.Wu \paper Any knot complement
covers at most one knot complement \jour Pacific J. Math. \vol 158
\pages 387--395 \yr 1993 \endref

\ref \key WWu2 \bysame \paper Covering invariant of graph
manifolds and cohopficity of 3-manifold groups \jour Proc. London
Math. Soc.\vol 68 \pages 221--242 \yr 1994 \endref

\ref \key WZh1 \by S.C.Wang and Q.Zhou \paper Embeddings of
Kleinian groups with torsion \jour Acta Math. Sinica \vol 17
\pages 21--34 \yr 2001 \endref

\ref \key WZh2 \bysame  \paper Any 3-manifold 1-dominates at most
finitely many geometric 3-manifolds \jour Math. Ann. \vol 332
\pages 525--535 \yr 2002
\endref

\ref \key YW \by F.Yu and S.C.Wang \paper Covering degrees are
uniquely determined by graph manifolds involved \jour Comm. Math.
Helv. \vol 74 \pages 238--247 \yr 1999 \endref

\enddocument